\documentclass[12pt]{article}

\usepackage{amsmath, epsfig, cite}
\usepackage{amssymb}
\usepackage{amsfonts}
\usepackage{latexsym}
\usepackage{float}

\newtheorem{thm}{Theorem}[section]

\newtheorem{lem}[thm]{Lemma}

\numberwithin{equation}{section}

\newcommand{\qed}{{\hfill$\square$}\medskip}

\begin{document}

%\linenumbers

\begin{center}
{\Large\bf Proof of some divisibility results on sums involving\\[5pt]
binomial coefficients}
\end{center}

\vskip 2mm \centerline{Ji-Cai Liu}
\begin{center}
{\footnotesize Department of Mathematics, East China Normal University, Shanghai 200241, PR China\\
{\tt jc2051@163.com} }
\end{center}

%%date: November 27, 2014
%%date : December 4, 2014

\vskip 0.7cm \noindent{\bf Abstract.}
By using the Rodriguez-Villegas-Mortenson supercongruences, we prove four supercongruences on sums involving binomial coefficients, which were originally conjectured by Sun.
We also confirm a related conjecture of Guo on integer-valued polynomials.

\vskip 3mm \noindent {\it Keywords}:
Supercongruences; Delannoy number; Legendre symbol; Zeilberger algorithm
\vskip 2mm
\noindent{\it MR Subject Classifications}: Primary 11A07; Secondary 33C05

\section{Introduction}
In 2003, Rodriguez-Villegas \cite{RV} conjectured 22 supercongruences for hypergeometric Calabi-Yau manifolds of dimension $d\le 3$. For manifolds of dimension $d=1$, associated to certain elliptic curves, four conjectural supercongruences were posed. Mortenson \cite{Mort1,Mort2} first proved these four supercongruences by using the Gross-Koblitz formula.

\begin{thm}\label{thm1} (Rodriguez-Villegas-Mortenson)
Suppose $p\ge 5$ is a prime. Then
\begin{align*}
\sum_{k=0}^{p-1}\frac{(1/2)_k^2}{(1)_k^2}\equiv \left(\frac{-1}{p}\right) \pmod{p^2},\quad
\sum_{k=0}^{p-1}\frac{(1/3)_k(2/3)_k}{(1)_k^2}\equiv \left(\frac{-3}{p}\right) \pmod{p^2},\\
\sum_{k=0}^{p-1}\frac{(1/4)_k(3/4)_k}{(1)_k^2}\equiv \left(\frac{-2}{p}\right) \pmod{p^2},\quad
\sum_{k=0}^{p-1}\frac{(1/6)_k(5/6)_k}{(1)_k^2}\equiv \left(\frac{-1}{p}\right) \pmod{p^2},
\end{align*}
where $\left(\frac{\cdot}{p}\right)$ denotes the Legendre symbol and $(x)_k=x(x+1)\cdots (x+k-1)$.
\end{thm}

Sun \cite{Sun} introduced the following two kinds of polynomials:
\begin{align*}
d_n(x)=\sum_{k=0}^n {n\choose k}{x\choose k}2^k \quad \text{and} \quad
s_n(x)=\sum_{k=0}^n {n\choose k}{x\choose k}{x+k\choose k}.
\end{align*}
Note that $d_n(m)$ are the Delannoy numbers, which count the number of paths from $(0,0)$ to $(m,n)$, only using steps $(1,0), (0,1)$ and $(1,1)$. For more information on Delannoy numbers, one can refer to \cite{BS}.

The first aim of this paper is to prove the following result, which was originally conjectured by Sun \cite[Conjecture 6.11]{Sun}.
\begin{thm} \label{conj2}
Suppose $p\ge 5$ is a prime. Then
\begin{align}
\sum_{k=0}^{p-1}(2k+1)s_k\left(-\frac{1}{2}\right)^2&\equiv \frac{3}{4}\left(\frac{-1}{p}\right)p^2 \pmod{p^4},\label{szw1}\\
\sum_{k=0}^{p-1}(2k+1)s_k\left(-\frac{1}{3}\right)^2&\equiv \frac{7}{9}\left(\frac{-3}{p}\right)p^2 \pmod{p^4},\label{szw2}\\
\sum_{k=0}^{p-1}(2k+1)s_k\left(-\frac{1}{4}\right)^2&\equiv \frac{13}{16}\left(\frac{-2}{p}\right)p^2 \pmod{p^4},\label{szw3}\\
\sum_{k=0}^{p-1}(2k+1)s_k\left(-\frac{1}{6}\right)^2&\equiv \frac{31}{36}\left(\frac{-1}{p}\right)p^2 \pmod{p^4}.\label{szw4}
\end{align}
\end{thm}

Recently, Guo \cite[Theorem 5.1]{Guo2} showed that for any odd prime $p$ and $p$-adic integer $x$,
\begin{align}
\sum_{k=0}^{p-1}(2k+1)s_k(x)^2\equiv p^2\sum_{k=0}^{p-1}\sum_{j=0}^k
\frac{(-1)^k}{k+1}{x+k\choose 2k}{x\choose j}{x+j\choose j}{2k\choose j+k}\pmod{p^4}.\label{bb1}
\end{align}

Recall that a polynomial $P(x)$ with real coefficients is called {\it integer-valued}, if $P(x)$ takes integer values for all $x\in \mathbb{Z}$.
The second aim of this paper is the prove the following result, which was originally conjectured by Guo \cite[Conjecture 5.5]{Guo2}.
\begin{thm}\label{conj3}
Let $n$ and $m$ be positive integers and $\varepsilon=\pm 1$. Then
\begin{align}
\frac{1}{n}\sum_{k=0}^{n-1}{\varepsilon}^k(2k+1)d_k(x)^m s_k(x)^m\label{aa1}
\end{align}
is integer-valued.
\end{thm}

\section{Proof of Theorem \ref{conj2}}
We need the following lemma (see \cite[Theorem 1.3]{Liu2}).
\begin{lem}\label{lem2}
Suppose $p\ge 5$ is a prime. Then
\begin{align*}
\sum_{k=0}^{2p-1}\frac{(1/2)_k^2}{(1)_k^2}\equiv \frac{5}{4}\left(\frac{-1}{p}\right) \pmod{p^2},\quad
\sum_{k=0}^{2p-1}\frac{(1/3)_k(2/3)_k}{(1)_k^2}\equiv \frac{11}{9}\left(\frac{-3}{p}\right) \pmod{p^2},\\
\sum_{k=0}^{2p-1}\frac{(1/4)_k(3/4)_k}{(1)_k^2}\equiv \frac{19}{16}\left(\frac{-2}{p}\right) \pmod{p^2},\quad
\sum_{k=0}^{2p-1}\frac{(1/6)_k(5/6)_k}{(1)_k^2}\equiv \frac{41}{36}\left(\frac{-1}{p}\right) \pmod{p^2}.
\end{align*}
\end{lem}

{\noindent\it Proof of Theorem \ref{conj2}.}
We begin with the following identity \cite[(2.5)]{Guo1}:
\begin{align}
{x\choose k}{x+k\choose k}{x\choose j}{x+j\choose j}=\sum_{s=0}^{j+k}{j+k\choose s}{s\choose j}{s\choose k}{x\choose s}{x+s\choose s}.\label{cc1}
\end{align}
Note that
\begin{align}
{x+k\choose 2k}{x\choose j}{x+j\choose j}={x\choose k}{x+k\choose k}{x\choose j}{x+j\choose j}\left.{\bigg/}\right.{2k\choose k}.\label{cc2}
\end{align}
Substituting \eqref{cc1} and \eqref{cc2} into the right-hand side of \eqref{bb1} and then exchanging the summation order gives
\begin{align}
&\sum_{k=0}^{p-1}(2k+1)s_k(x)^2\notag\\
&\equiv p^2\sum_{s=0}^{2p-2}\sum_{k=0}^{p-1}\sum_{j=0}^k
\frac{(-1)^k}{k+1}{2k\choose j+k}{j+k\choose s}{s\choose j}{s\choose k}{x\choose s}{x+s\choose s}{\bigg/}{2k\choose k} \pmod{p^4}.\label{cc3}
\end{align}
By the Chu-Vandermonde identity, we have
\begin{align}
\sum_{j=0}^k{2k\choose j+k}{j+k\choose s}{s\choose j}=\sum_{j=0}^k{2k\choose s}{s\choose j}{2k-s\choose k-j}
={2k\choose s}{2k\choose k}.\label{cc4}
\end{align}
It follows from \eqref{cc3} and \eqref{cc4} that
\begin{align}
\sum_{k=0}^{p-1}(2k+1)s_k(x)^2
\equiv p^2\sum_{s=0}^{2p-2}\sum_{k=0}^{p-1}
\frac{(-1)^k}{k+1}{2k\choose s}{s\choose k}{x\choose s}{x+s\choose s} \pmod{p^4}.\label{cc5}
\end{align}

Applying the following identity \cite[(2.6)]{Liu}:
\begin{align*}
\sum_{k=0}^{s}\frac{(-1)^k}{k+1}{2k\choose s}{s\choose k}=(-1)^s,
\end{align*}
we obtain
\begin{align}
p^2\sum_{s=0}^{p-1}\sum_{k=0}^{p-1}
\frac{(-1)^k}{k+1}{2k\choose s}{s\choose k}{x\choose s}{x+s\choose s}
=p^2\sum_{s=0}^{p-1}(-1)^s{x\choose s}{x+s\choose s}.\label{cc6}
\end{align}

On the other hand, we have the following supercongruence \cite[(3.2)]{Liu}:
\begin{align}
\sum_{k=0}^{p-1}
\frac{(-1)^k}{k+1}{2k\choose s}{s\choose k}\equiv (-1)^s\left(-1+\frac{2p}{s+1}\right) \pmod{p^2}\label{cc7}
\end{align}
for $p\le s\le 2p-2$. It follows that
\begin{align}
&p^2\sum_{s=p}^{2p-2}\sum_{k=0}^{p-1}
\frac{(-1)^k}{k+1}{2k\choose s}{s\choose k}{x\choose s}{x+s\choose s}\notag\\
&\equiv 2p^3\sum_{s=p}^{2p-1}\frac{(-1)^s}{s+1}{x\choose s}{x+s\choose s}- p^2\sum_{s=p}^{2p-1}(-1)^s{x\choose s}{x+s\choose s} \pmod{p^4},\label{cc8}
\end{align}
where we have used the fact that ${x\choose s}{x+s\choose s}\equiv 0\pmod{p^2}$ for $s=2p-1$ and
$x=-\frac{1}{2},-\frac{1}{3},-\frac{1}{4},-\frac{1}{6}$.

Note that
\begin{align*}
\sum_{s=0}^{n-1}\frac{(-1)^s}{s+1}{x\choose s}{x+s\choose s}=\frac{n(-1)^{n+1}}{x(x+1)}{x\choose n}{x+n\choose n},
\end{align*}
which can be easily proved by induction on $n$.
So we have
\begin{align}
\sum_{s=p}^{2p-1}\frac{(-1)^s}{s+1}{x\choose s}{x+s\choose s}&=\sum_{s=0}^{2p-1}\frac{(-1)^s}{s+1}{x\choose s}{x+s\choose s}-\sum_{s=0}^{p-1}\frac{(-1)^s}{s+1}{x\choose s}{x+s\choose s}\notag\\
&\equiv 0 \pmod{p}\label{cc9}
\end{align}
for $x=-\frac{1}{2},-\frac{1}{3},-\frac{1}{4},-\frac{1}{6}$.

Furthermore, combining \eqref{cc5}, \eqref{cc6}, \eqref{cc8} and \eqref{cc9}, we get
\begin{align}
&\sum_{k=0}^{p-1}(2k+1)s_k(x)^2\notag\\
&\equiv p^2\left(2\sum_{s=0}^{p-1}(-1)^s{x\choose s}{x+s\choose s}-\sum_{s=0}^{2p-1}(-1)^s{x\choose s}{x+s\choose s}\right) \pmod{p^4}\label{cc10}
\end{align}
for $x=-\frac{1}{2},-\frac{1}{3},-\frac{1}{4},-\frac{1}{6}$.

Finally, noting that
\begin{align*}
(-1)^s{x\choose s}{x+s\choose s}=\frac{(-x)_s (1+x)_s}{(1)_s^2},
\end{align*}
and then using Theorem \ref{thm1} and Lemma \ref{lem2}, we obtain \eqref{szw1}-\eqref{szw4}.
This completes the proof.
\qed

\section{Proof of Theorem \ref{conj3}}
Chen and Guo \cite{CG} defined the following multi-variable Schmidt polynomials:
\begin{align*}
S_n(x_0,\cdots,x_n)=\sum_{k=0}^n {n+k\choose 2k}{2k\choose k}x_k.
\end{align*}
The following lemma is a special case of \cite[Theorem 1.1]{CG}. It has already been
used by Guo \cite{Guo2} to prove Sun's conjectures on integer-valued polynomials.
\begin{lem}\label{lem3}
Let $n$ and $m$ be positive integers. Then all the coefficients in
\begin{align*}
\sum_{k=0}^{n-1}{\varepsilon}^k(2k+1)S_k(x_0,\cdots,x_k)^m
\end{align*}
are multiples of $n$.
\end{lem}

\begin{lem}
For any non-negative integer $n$, we have
\begin{align}
d_n(x)s_n(x)=\sum_{k=0}^n {n+k\choose 2k}{2k\choose k}\sum_{j=0}^k \sum_{i=0}^j
{x+j\choose k+j}{x\choose i}{k\choose j}{j\choose i}2^i. \label{bb2}
\end{align}
\end{lem}
\noindent{\it Proof.}
Since both sides of \eqref{bb2} are polynomials in $x$ of degree $3n$, it suffices to prove that for any non-negative integers $n$ and $m$,
\begin{align}
&\sum_{i=0}^n \sum_{j=0}^n {n\choose i}{m\choose i}{n\choose j}{m\choose j}{m+j\choose j}2^i\notag\\
&=\sum_{k=0}^n\sum_{j=0}^n \sum_{i=0}^n {n+k\choose 2k}{2k\choose k}{m+j\choose k+j}{m\choose i}{k\choose j}{j\choose i}2^i. \label{bb3}
\end{align}
Clearly, \eqref{bb3} is equivalent to
\begin{align}
&\sum_{i=0}^m \sum_{j=0}^m {n\choose i}{m\choose i}{n\choose j}{m\choose j}{m+j\choose j}2^i\notag\\
&=\sum_{k=0}^m\sum_{j=0}^m \sum_{i=0}^m {n+k\choose 2k}{2k\choose k}{m+j\choose k+j}{m\choose i}{k\choose j}{j\choose i}2^i. \label{bb4}
\end{align}
Let $A_m^{(1)}(n)$ and $A_m^{(2)}(n)$ denote the left-hand side and the right-hand side of \eqref{bb4}, respectively. Applying the multi-Zeilberger algorithm \cite{AZ,PWZ}, we find that $A_m^{(1)}(n)$ and $A_m^{(2)}(n)$ satisfy the same recurrence of order $4$:
\begin{align*}
&(m+1)^3(m+2)(3m^2+18m+26)A_m^{(s)}(n)-2(m+2)(12m^3n^2+12m^3n+90m^2n^2\\
&+3m^3+90m^2n+212mn^2+23m^2+212mn+156n^2+55m+156n+41)A_{m+1}^{(s)}(n)\\
&-2(3m^6+30m^4n^2+45m^5+30m^4n+300m^3n^2+287m^4+300m^3n+1094m^2n^2\\
&+995m^3+1094m^2n+1720mn^2+1964m^2+1720mn+978n^2+2070m+978n\\
&+898)A_{m+2}^{(s)}(n)-2(m+3)(12m^3n^2+12m^3n+90m^2n^2+3m^3+90m^2n+212mn^2\\
&+22m^2+212mn+154n^2+50m+154n+34)A_{m+3}^{(s)}(n)+(m+3)(m+4)^3\\
&\times(3m^2+12m+11)A_{m+4}^{(s)}(n)=0,\quad \text{for $s=1,2.$}
\end{align*}
It is easily checked that $A_{m}^{(1)}(n)=A_{m}^{(2)}(n)$ for $0\le m\le 3$.
This proves \eqref{bb4}. \qed

{\noindent\it Proof of Theorem \ref{conj3}.}
We can rewrite \eqref{bb2} as
\begin{align}
d_n(x)s_n(x)=\sum_{k=0}^n {n+k\choose 2k}{2k\choose k}f_k(x),\label{dd1}
\end{align}
where
\begin{align*}
f_k(x)=\sum_{j=0}^k \sum_{i=0}^j{x+j\choose k+j}{x\choose i}{k\choose j}{j\choose i}2^i.
\end{align*}
Clearly, these $f_k(x)$ are integer-valued polynomials for $0\le k\le n$.
From Lemma \ref{lem3} and the identity \eqref{dd1}, we immediately conclude that $\eqref{aa1}$ is integer-valued.
\qed


\begin{thebibliography}{99}
\small \setlength{\itemsep}{-.8mm}

\bibitem{AZ}M. Apagodu and D. Zeilberger, Multi-variable Zeilberger and Almkvist--Zeilberger
algorithms and the sharpening of Wilf--Zeilberger theory, Adv. Appl. Math. 37 (2006), 139--152.

\bibitem{BS}C. Banderier and S.R. Schwer, Why Delannoy's numbers? J. Stat. Plan. Infer. 135 (2005), 40--54.

\bibitem{CG}Q.-F. Chen and V.J.W. Guo, On the divisibility of sums involving powers of multi-variable Schmidt polynomials, preprint, 2014, arXiv:1412.5734.

\bibitem{Guo1}V.J.W. Guo, Some congruences involving powers of Legendre polynomials,
Integral Transforms Spec. Funct. 26 (2015), 660--666.

\bibitem{Guo2}V.J.W. Guo, Proof of Sun's conjectures on integer-valued polynomials,
 J. Math. Anal. Appl. 444 (2016), 182--191.

\bibitem{Liu}J.-C. Liu, A supercongruence involving Delannoy numbers and Schr\"oder numbers, J. Number Theory 168 (2016), 117--127.

\bibitem{Liu2}J.-C. Liu, Congruences for truncated hypergeometric series ${}_2F_1$, Bull. Aust. Math. Soc., online, doi: 10.1017/S0004972717000181.

\bibitem{Mort1}E. Mortenson, A supercongruence conjecture of Rodriguez-Villegas for a certain truncated hypergeometric function, J. Number Theory 99 (2003), 139--147.

\bibitem{Mort2}E. Mortenson, Supercongruences between truncated ${}_2F_1$ hypergeometric functions and their Gaussian analogs, Trans. Amer. Math. Soc. 355 (2003), 987--1007.

\bibitem{PWZ}M. Petkov\v{s}ek, H. S. Wilf and D. Zeilberger, $A=B$, A K Peters, Ltd., Wellesley, MA, 1996.

\bibitem{RV} F. Rodriguez-Villegas, Hypergeometric families of Calabi-Yau manifolds, in: Calabi-Yau Varieties and Mirror Symmetry (Toronto, ON, 2001), Fields Inst. Commun., 38, Amer.
Math. Soc., Providence, RI, 2003, 223--231.

\bibitem{Sun}Z.-W. Sun, Supercongruences involving dual sequences, preprint, 2015, arXiv:1512.00712.

\end{thebibliography}
\end{document}